\documentclass[a4paper,11pt]{article}

\usepackage{amsmath}
\usepackage{amsfonts}
\usepackage{amsthm}
\usepackage[english]{babel}

\newtheorem{lemma}{Lemma}[section]
\newtheorem{definition}{Definition}[section]
\newtheorem{theorem}{Theorem}[section]
\newtheorem{proposition}{Proposition}[section]
\newtheorem{corollary}{Corollary}[section]
\newtheorem{remark}{Remark}[section]

\newcommand{\mR}{\mathbb{R}}
\newcommand{\mC}{\mathbb{C}}

\newcommand{\mE}{\mathbb{E}}

\newcommand{\mP}{\mathbb{P}}

\newcommand{\cP}{\mathcal{P}}

\newcommand{\ux}{\underline{x}}

\begin{document}

\title{Clifford-Gegenbauer polynomials related to the Dunkl Dirac operator}

\author{H.\ De Bie\thanks{Postdoctoral Fellow of the Research Foundation - Flanders (FWO), E-mail: {\tt Hendrik.DeBie@UGent.be}}, \and N.\ De Schepper\footnote{E-mail: {\tt nds@cage.ugent.be}}}

\date{\small{Clifford Research Group -- Department of Mathematical Analysis}\\
\small{Faculty of Engineering -- Ghent University\\ Krijgslaan 281, 9000 Gent,
Belgium}}

\maketitle

\begin{abstract}
We introduce the so-called Clifford-Gegenbauer polynomials in the framework of Dunkl operators, as well on the unit ball $B(1)$, as on the Euclidean space $\mathbb{R}^m$. In both cases we obtain several properties of these polynomials, such as a Rodrigues formula, a differential equation and an explicit relation connecting them with the Jacobi polynomials on the real line. As in the classical Clifford case, the orthogonality of the polynomials on $\mathbb{R}^m$ must be treated in a completely different way than the orthogonality of their counterparts on $B(1)$. In case of $\mathbb{R}^m$, it must be expressed in terms of a bilinear form instead of an integral. Furthermore, in this paper the theory of Dunkl monogenics is further developed.
\end{abstract}



\section{Introduction}
Dunkl operators (see \cite{MR951883, Dunkl}) are combinations of differential and difference operators, associated to a finite reflection group $G$. One of the interesting aspects of these operators is that they allow for the construction of a Dunkl Laplacian, which is a combination of the classical Laplacian in $\mathbb{R}^m$ with some difference terms, such that the resulting operator is only invariant under $G$ and not under the whole orthogonal group. Moreover, they are directly related to quantum integrable models of Calogero type (see e.g. \cite{Ols}) and have as such received a lot of attention in the physics literature.

Clifford analysis (see \cite{rood,groen}), in its most basic form, is a refinement of the theory of harmonic analysis in $m$-dimensional Euclidean space. By introducing the so-called Dirac operator, the square of which equals the Laplace operator, one introduces the notion of monogenic functions. These are, at the same time, a refinement of harmonic functions and a generalization of holomorphic functions in one complex variable.\\
Generalizations of the classical Gegenbauer polynomials to the Clifford analysis framework are called Clifford-Gegenbauer polynomials and were introduced as well on the closed unit ball $B(1)$ (see \cite{Cnops}), as on the Euclidean space $\mathbb{R}^m$ (see \cite{wave,WCAA7}). For the superspace case, which can be seen as the study of differential operators invariant under the action of the group $O(m) \times \mathrm{Sp}(2n)$, we refer to \cite{DBS3}.\\
In this paper, we adapt the definition of the Clifford-Gegenbauer polynomials, both on $B(1)$, as on $\mathbb{R}^m$, to the case of Dunkl operators.

The paper is organized as follows. In Section 2 we first give some background on  Dunkl operators. Then we prove some fundamental results concerning the Dunkl Dirac operator and its nullsolutions, called Dunkl monogenics. In Section 3 we introduce Clifford-Gegenbauer polynomials on $B(1)$ related to the Dunkl Dirac operator. Basic properties, such as a Rodrigues formula, a differential equation, recurrence relations and an orthogonality relation are derived. Moreover, we obtain an expression of these newly introduced polynomials in terms of the Jacobi polynomials on the real line. Next, the Clifford-Gegenbauer polynomials on $\mathbb{R}^m$ are adapted to the case of Dunkl operators (Section 4). They satisfy similar properties as their counterparts on the unit ball. However, the orthogonality must be treated completely different; it is expressed in terms of a bilinear form instead of an integral.
\section{Clifford Dunkl setting}

\subsection{Dunkl operators}

Denote by $\langle .,. \rangle$ the standard Euclidean scalar product in $\mR^{m}$ and by $|x| = \langle x, x\rangle^{1/2}$ the associated norm. For $\alpha \in \mR^{m} - \{ 0\}$, the reflection $r_{\alpha}$ in the hyperplane orthogonal to $\alpha$ is given by
\[
r_{\alpha}(x) = x - 2 \frac{\langle \alpha, x\rangle}{|\alpha|^{2}}\alpha, \quad x \in \mR^{m}.
\]

A root system is a finite subset $R \subset \mR^{m}$ of non-zero vectors such that, for every $\alpha \in R$, the associated reflection $r_{\alpha}$ preserves $R$. We will assume that $R$ is reduced, i.e. $R \cap \mR \alpha = \{ \pm \alpha\}$ for all $\alpha \in R$. Each root system can be written as a disjoint union $R = R_{+} \cup (-R_{+})$, where $R_{+}$ and $-R_{+}$ are separated by a hyperplane through the origin. The subgroup $G \subset O(m)$ generated by the reflections $\{r_{\alpha} | \alpha \in R\}$ is called the finite reflection group associated with $R$. We will also assume that $R$ is normalized such that $\langle \alpha, \alpha\rangle = 2$ for all $\alpha \in R$. For more information on finite reflection groups we refer the reader to \cite{Humph}.

A multiplicity function $k$ on the root system $R$ is a $G$-invariant function $k: R \rightarrow \mC$, i.e. $k(\alpha) = k(h \alpha)$ for all $h \in G$. We will denote $k(\alpha)$ by $k_{\alpha}$.

Fixing a positive subsystem $R_{+}$ of the root system $R$ and a multiplicity function $k$, we introduce the Dunkl operators $T_{i}$ associated to $R_{+}$ and $k$ by (see \cite{MR951883, Dunkl})
\[
T_{i} f(x)= \partial_{x_{i}} f(x) + \sum_{\alpha \in R_{+}} k_{\alpha} \alpha_{i} \frac{f(x) - f(r_{\alpha}(x))}{\langle \alpha, x\rangle}, \qquad f \in C^{1}(\mR^{m}).
\]
An important property of the Dunkl operators is that they commute, i.e. $T_{i} T_{j} = T_{j} T_{i}$.

The Dunkl Laplacian is given by $\Delta_{k} = \sum_{i=1}^{m} T_i^2$, or more explicitly by
\[
\Delta_{k} f(x) = \Delta f(x) + 2 \sum_{\alpha \in R_{+}} k_{\alpha} \left( \frac{\langle \nabla f(x), \alpha \rangle}{\langle \alpha, x \rangle}  - \frac{f(x) - f(r_{\alpha}(x))}{\langle \alpha, x \rangle^{2}} \right)
\]
with $\Delta$ the classical Laplacian and $\nabla$ the gradient operator.


If we let $\Delta_{k}$ act on $|x|^2$ we find $\Delta_{k} \lbrack |x|^2 \rbrack = 2m + 4 \gamma = 2 \mu$, where $\gamma = \sum_{\alpha \in R_+} k_{\alpha}$. We call $\mu$ the Dunkl dimension, because most special functions related to $\Delta_{k}$ behave as if one would be working with the classical Laplace operator in a space with dimension $\mu$. 


The operators
\[ E:= \frac{1}{2}|x|^2,\>\> F:= -\frac{1}{2}\Delta_{k} \quad\text{and}\>\> 
H:= \mE +\mu/2
\]
where $\mathbb{E} := \sum_{i=1}^m x_i \partial_{x_i}$ is the Euler operator, satisfy the defining relations of the Lie algebra $\mathfrak{sl}_2$ (see e.g. \cite{He}). 
They are given by
\begin{equation}
\label{sl2reldunkl}
\big[H,E\big] = 2E,\>\> \big[H,F\big] = -2F,\>\>\big[E,F\big] = H.
\end{equation}

We also have the following important property of the Dunkl operators
\begin{equation}\label{prod}
T_i \lbrack f g \rbrack = T_i \lbrack f \rbrack \ g + f \ T_i \lbrack g \rbrack \ \ , \ \ i=1,2,\ldots,m
\end{equation}
if $f$ or $g$ are $G$-invariant.


\subsection{Dunkl Dirac operators}

For the sequel of this paper we restrict ourselves to multiplicity functions satisfying $k_{\alpha} \geq 0$, $\forall \alpha \in R_+$; hence $\gamma = \sum_{\alpha \in R_+} k_{\alpha}  \geq 0$ and for the Dunkl dimension $\mu$ we have that $\mu = m + 2 \gamma >1$.

From now on we consider functions $f: \mathbb{R}^m \rightarrow \mathbb{R}_{0,m}$. Hereby, $\mathbb{R}_{0,m}$ denotes the Clifford algebra over $\mR^{m}$ generated by $e_{i}, i=1, \ldots,m$, subject to the relations $e_{i} e_{j} + e_{j} e_{i} = -2 \delta_{ij}$.

In what follows, the bar denotes the Clifford conjugation; an anti-involution for which $\overline{e_i} = -e_i$ ($i=1,\ldots,m$).

$D_k$ stands for the Dunkl Dirac operator $D_k = \sum_{j=1}^m e_j T_j$ with $T_j$ the Dunkl operators. $\ux = \sum_{j=1}^m e_j x_j$ is the so-called vector variable. It is clear that $D_{k}^{2} = - \Delta_{k}$ and that $\ux^{2} = -|\ux|^{2}= -r^{2}$.

By definition we have that for $f(\ux) = \sum_{A} e_A f_A(\ux)$ with $f_A : \mathbb{R}^m \rightarrow \mathbb{R}$ :
\begin{displaymath}
D_k \lbrack f(\ux) \rbrack = \sum_{j=1}^m \sum_A e_j e_A T_j \lbrack f_A \rbrack \qquad \mathrm{and} \qquad \lbrack f(\ux) \rbrack D_k = \sum_{j=1}^m \sum_A T_j \lbrack f_A \rbrack e_A e_j  .
\end{displaymath}
A function $f$ which satisfies $D_k \lbrack f \rbrack = 0$ is called left Dunkl monogenic, while a function $f$ which satisfies $\lbrack f \rbrack D_k = 0$ is called right Dunkl monogenic.\\
A left Dunkl monogenic homogeneous polynomial $M_k$ of degree $k$ ($k \geq 0$) in $\mathbb{R}^m$ is called a left solid inner Dunkl monogenic of order $k$. A left Dunkl monogenic homogeneous function $Q_k$ of degree $-(k+\mu-1)$ in $\mathbb{R}^m \setminus \lbrace 0 \rbrace$ is called a left solid outer Dunkl monogenic of order $k$. The restriction to $S^{m-1}$ of a left solid inner Dunkl monogenic is called a left inner Dunkl monogenic, while the restriction to $S^{m-1}$ of a left solid outer Dunkl monogenic is called a left outer Dunkl monogenic. The set of all left solid inner Dunkl monogenics of order $k$ will be denoted by $M_{\ell}^+(k)$, while the set of all left solid outer Dunkl monogenics of order $k$ will be denoted by $M_{\ell}^-(k)$. Moreover, the space of left inner Dunkl monogenics of order $k$ is denoted by ${\mathcal M}_{\ell}^+(k)$, while the space of left outer Dunkl monogenics is denoted by ${\mathcal M}_{\ell}^-(k)$. Similar definitions hold in case of right Dunkl monogenics.

The angular Dunkl Dirac operator (or Gamma operator) is defined as $\Gamma_k := D_k \ux + \mu + \mathbb{E}$.\\
Taking into account the relation (see for e.g. \cite{Hendrik, Vladimir})
\begin{equation}\label{fundamenteel}
\lbrace D_k , \ux \rbrace = D_k \ux + \ux D_k = - ( 2 \mathbb{E} + \mu)  ,
\end{equation}
we easily obtain the following proposition.
\begin{proposition}\label{decxdirac}~\\
One has that
\begin{displaymath}
\ux D_k = - \mathbb{E}-\Gamma_k \qquad \mathrm{or} \qquad \Gamma_k = -\ux D_k - \mathbb{E} .
\end{displaymath}
\end{proposition}
\textit{Proof.} We have consecutively
\begin{displaymath}
\ux D_k  =  - D_k \ux - 2 \mathbb{E}-\mu  =  - \Gamma_k + \mu+\mathbb{E}-2 \mathbb{E} - \mu =  - \Gamma_k - \mathbb{E}  . \ \ \ \ \ \square
\end{displaymath}
Now we immediately obtain that $\mathbb{E}$ measures the degree of homogeneity, while $\Gamma_k$ measures the degree of Dunkl monogenicity.
\begin{proposition}\label{measure}~\\
(i) $R_k \in {\mathcal P}_k$ : $\mathbb{E} \lbrack R_k \rbrack = k R_k$\\
(ii) $M_k \in M_{\ell}^+(k)$ : $ \Gamma_k \lbrack M_k \rbrack = -k M_k$\\
(iii) $M_k \in M_{\ell}^+(k)$ : $ \Gamma_k \lbrack \ux \ M_k \rbrack = (k+\mu-1) \ \ux  M_k$.
\end{proposition}
\textit{Proof.}\\
(i) Straightforward.\\
(ii) $\Gamma_k \lbrack M_k \rbrack = (- \mathbb{E} - \ux D_k) \lbrack M_k \rbrack = -k M_k$.\\
(iii) By means of (i) and (ii) we obtain
\begin{displaymath}
\ux D_k \lbrack \ux M_k \rbrack = \ux \ \left( \Gamma_k - \mu - \mathbb{E} \right) \lbrack M_k \rbrack = \ux \ (-k-\mu-k) \ M_k = (-2k-\mu) \ \ux M_k  .
\end{displaymath}
Using Proposition \ref{decxdirac} and the fact that $\ux M_k \in {\mathcal P}_{k+1}$, gives
\begin{displaymath}
(-2k-\mu) \ \ux M_k = (- \mathbb{E}-\Gamma_k) \lbrack \ux M_k \rbrack = -(k+1) \ \ux M_k - \Gamma_k \lbrack \ux M_k \rbrack  ,
\end{displaymath}
which yields the desired result. \ \ \ \ \ \ $\square$
\\
The following lemma is proved by induction using (\ref{fundamenteel}).
\begin{lemma}\label{belangrijk}~\\
For $M_k \in M_{\ell}^+(k)$ we have that
\begin{equation*}
D_k \lbrack \ux^s M_k \rbrack = 
  \begin{cases}
       -s \ux^{s-1} M_k & \text{if $s$ is even},\\
       -(s-1+2k+\mu) \ux^{s-1} M_k & \text{if $s$ is odd}.
  \end{cases}
\end{equation*}
\end{lemma}
Now the following connection between left solid inner and outer Dunkl monogenics can be proved.
\begin{proposition}\label{connectionmonogenic}~\\
(i) For $M_k \in M_{\ell}^+(k)$ we have that
\begin{displaymath}
Q_k(\ux) = \frac{\ux}{|\ux|^{\mu}} \ M_k \left( \frac{\ux}{|\ux|^2} \right) = \frac{\ux}{|\ux|^{\mu+2k}} \ M_k(\ux) \in M_{\ell}^-(k)  .
\end{displaymath}
(ii) For $\widetilde{Q_k} \in M_{\ell}^-(k)$ we have that
\begin{displaymath}
\widetilde{M_k}(\ux) = \frac{\ux}{|\ux|^{\mu}} \ \widetilde{Q_k} \left( \frac{\ux}{|\ux|^2} \right) = \ux \ |\ux|^{2k+\mu-2} \ \widetilde{Q_k}(\ux) \in M_{\ell}^+(k)  .
\end{displaymath}
Hence on the unit sphere $S^{m-1}$ we have the following connection
\begin{displaymath}
M_k \in {\mathcal M}_{\ell}^+(k) \ \Longleftrightarrow \ \underline{\omega} \ M_k(\underline{\omega}) \in {\mathcal M}_{\ell}^-(k) \quad , \quad \underline{\omega} \in S^{m-1} .
\end{displaymath}
\end{proposition}
\textit{Proof.}\\
(i) Clearly $Q_k$ is homogeneous of degree $-(k+\mu-1)$. Moreover taking into account the product rule (\ref{prod}), we obtain that in $\mathbb{R}^m \setminus \lbrace 0 \rbrace$ : 
\begin{eqnarray*}
D_k \lbrack Q_k(\ux) \rbrack & = & D_k \left\lbrack \frac{1}{|\ux|^{\mu+2k}} \right\rbrack \ \ux M_k + \frac{1}{|\ux|^{\mu+2k}} \ D_k \lbrack \ux M_k \rbrack\\
& = & \sum_{i=1}^m e_i \partial_{x_i} \left\lbrack \frac{1}{|\ux|^{\mu+2k}} \right\rbrack  \ \ux M_k - \frac{(\mu+2k)}{|\ux|^{\mu+2k}} \ M_k\\
& = & - \sum_{i=1}^m e_i \ (\mu+2k) \ \frac{1}{|\ux|^{\mu+2k+1}} \frac{x_i}{|\ux|} \ \ux M_k -\frac{(\mu+2k)}{|\ux|^{\mu+2k}} \ M_k \ = \ 0 .\\
\end{eqnarray*}
(ii) The converse result is proved in a similar way. \ \ \ \ $\square$\vspace{0,2cm}\\
A straightforward calculation yields the following lemma.
\begin{lemma}\label{diracopx}~\\
One has $D_k \lbrack \ux \rbrack = -\mu$.
\end{lemma}
\textit{Proof.} From (\ref{fundamenteel}) we obtain that $\lbrace D_k , \ux \rbrace \lbrack 1 \rbrack = - (2 \mathbb{E}+\mu) \lbrack 1 \rbrack$ and hence $D_k \lbrack \ux \rbrack = - \mu$. \ \ \  \ $\square$\vspace{0,2cm}\\
The angular Dunkl Dirac operator only acts on the angular co-ordinates, whence its name.
\begin{proposition}\label{angular}~\\
One has that $\lbrack \Gamma_k , f(r) \rbrack = 0$.
\end{proposition}
\textit{Proof.} Let $G(\ux)$ denote a Clifford algebra-valued function. Taking into account Proposition \ref{decxdirac} and the product rule (\ref{prod}), we arrive at
\begin{eqnarray*}
& & - \Gamma_k \lbrack f(r) \ G(\ux) \rbrack\\
& = & (\mathbb{E} + \ux D_k) \lbrack f(r) \ G(\ux) \rbrack \ = \ \mathbb{E} \lbrack f(r) \rbrack \ G(\ux) + f(r) \ \mathbb{E} \lbrack G(\ux) \rbrack + \ux \ D_k \lbrack f(r) \ G(\ux) \rbrack\\
& = & r \partial_r \lbrack f(r) \rbrack \ G(\ux) + f(r) \ \mathbb{E} \lbrack G(\ux) \rbrack + \ux \ \partial_{\ux} \lbrack f(r) \rbrack \ G(\ux) + f(r) \ \ux \ D_k \lbrack G(\ux) \rbrack\\
& = & - f(r) \ \Gamma_k \lbrack G(\ux) \rbrack + r \partial_r \lbrack f(r) \rbrack \ G(\ux) + \ux \ \frac{\ux}{r} \ \partial_r \lbrack f(r) \rbrack \ G(\ux) \ = \ -f(r) \ \Gamma_k \lbrack G(\ux) \rbrack   ,
\end{eqnarray*}
which proves the statement. \ \ \ $\square$

We also have the following important decomposition (see \cite{Vladimir}).
\begin{theorem}[Fischer decomposition]
The space $\cP_{k}$ of homogeneous polynomials of degree $k$ taking values in $\mR_{0,m}$ decomposes as
\[
\cP_{k} = \bigoplus_{i=0}^{k} \ux^{i} M_{\ell}^+(k-i),
\]
where $M_{\ell}^+(k-i)$ is the space of left solid inner Dunkl monogenics of order $k-i$.
\end{theorem}

Using the Gamma operator, it is now possible to construct projection operators on each piece in the Fischer decomposition in a much easier way than was obtained in \cite{Vladimir}. Indeed, it is easy to see that each summand in this decomposition corresponds to a different eigenvalue of $\Gamma_{k}$. Hence, we can use the Gamma operator to construct projection operators on each summand of the decomposition. So, the operator $\mP^k_i$ (acting on $\cP_k$) defined by
\begin{eqnarray*}
\mP^k_i&=& \prod_{r=0, r \neq i/2}^{\lfloor \frac{k}{2} \rfloor} \dfrac{\Gamma_{k} +k -2r}{i-2r} \prod_{s=0}^{\lfloor \frac{k-1}{2} \rfloor} \dfrac{\Gamma_{k} -k + 2s +2 -\mu}{-2k + i + 2s +2 -\mu}, \quad \mbox{$i$ even}\\
&=&\prod_{r=0}^{\lfloor \frac{k}{2} \rfloor} \dfrac{\Gamma_{k} +k -2r}{2k-i+\mu-1-2r} \prod_{s=0, s \neq \frac{i-1}{2}}^{\lfloor \frac{k-1}{2} \rfloor} \dfrac{\Gamma_{k} -k + 2s +2 -\mu}{2s +1 - i}, \quad \mbox{$i$ odd}\\
\end{eqnarray*}
clearly satisfies
\[
\mP^k_i \left( \ux^j M_{\ell}^+(k-j) \right) = \delta_{ij} \ux^j M_{\ell}^+(k-j).
\]

\begin{proposition}\label{exprDirac}~\\
In terms of spherical co-ordinates the Dunkl Dirac operator takes the form
\begin{displaymath}
D_k = \underline{\omega} \left( \partial_r + \frac{1}{r} \Gamma_k \right) .
\end{displaymath}
\end{proposition}
\textit{Proof.} Again by means of Proposition \ref{decxdirac}, we obtain
\begin{displaymath}
r \underline{\omega} \ D_k = -( \mathbb{E}  +\Gamma_k) = -(r \partial_r +\Gamma_k) \qquad \mathrm{or} \qquad D_k = \underline{\omega} \left(\partial_r + \frac{1}{r} \Gamma_k \right)  . \ \ \ \ \square
\end{displaymath}
\begin{lemma}~\\
One has that $\Gamma_k \lbrack \underline{\omega} \rbrack  = (\mu-1) \underline{\omega}$,  $\underline{\omega} \in S^{m-1}$.
\end{lemma}
\textit{Proof.} Taking into account Proposition \ref{angular} and Proposition \ref{measure} (iii), we find consecutively
\begin{displaymath}
\Gamma_k \lbrack \underline{\omega} \rbrack  =  \Gamma_k \left\lbrack \frac{\ux}{r} \right\rbrack =  \frac{1}{r} \ \Gamma_k \lbrack \ux \rbrack =  \frac{1}{r} \ (\mu - 1) \ \ux =  (\mu-1) \ \underline{\omega} . \ \ \ \ \ \square
\end{displaymath}

In the sequel the following positive weight function will play a crucial role:
\begin{displaymath}
w_k(\ux) = \prod_{\alpha \in R_+} |\langle x,\alpha \rangle|^{2 k_{\alpha}} .
\end{displaymath}
It is homogeneous of degree $2 \gamma = 2 \ \sum_{\alpha \in R_+}  k_{\alpha}$ and invariant under reflections from the root system $R_+$.\\

A basic integral formula is the Stokes formula (see \cite{Paula}), the proof of which heavily relies on the $G$-invariance of the weight function $w_k(\ux)$.\\
Let $\Omega$ be a sufficiently smooth domain with boundary $\Gamma = \partial \Omega$ and define the oriented surface element $d \sigma (\ux)$ on $\Gamma$ by the Clifford differential form:
\begin{displaymath}
d\sigma(\ux) = \sum_{j=1}^m (-1)^j \ e_j \ dx^{M \setminus \lbrace j \rbrace} ,
\end{displaymath}
where 
\begin{displaymath}
dx^{M \setminus \lbrace j \rbrace} = dx_1 \wedge \ldots \wedge \lbrack dx_j \rbrack \wedge \ldots \wedge dx_m \quad , \quad j=1,2,\ldots,m  .
\end{displaymath}
If $n(\ux)$ stands for the outward pointing unit normal at $\ux \in \Gamma$, then $d\sigma(\ux) = n(\ux) \ d\Sigma(\ux)$, $d\Sigma(\ux)$ being the elementary Lebesgue surface measure.
\begin{theorem}[Stokes and Cauchy theorem]\label{stokes}~\\
Let $\Omega$ be a sufficiently smooth domain invariant under the action of $G$, $\Gamma = \partial \Omega$ and $f,g \in C^{\infty}(\Omega)$. Then
\begin{displaymath}
\int_{\Omega} \lbrack (f D_k) \ g + f \ (D_k g) \rbrack  \ w_k(\ux) \ dV(\ux) = \int_{\Gamma} f \ w_k(\ux) \ d\sigma(\ux) \ g  .
\end{displaymath}
Moreover, if $f$ is right monogenic in $\Omega$ and $g$ is left monogenic in $\Omega$, one has
\begin{displaymath}
\int_{\Gamma} f \ w_k(\ux) \ d\sigma(\ux) \ g = 0  .
\end{displaymath}
\end{theorem}
We now have all the necessary results at our disposal in order to prove the orthogonality of the inner and outer Dunkl monogenics.
\begin{theorem}[Orthogonality of Dunkl monogenics]\label{orthomon}~\\
(i) Left inner Dunkl monogenics of different degree are orthogonal, i.e. for $M_t \in {\mathcal M}^+_{\ell}(t)$ and $M_k \in {\mathcal M}^+_{\ell}(k)$ with $t \not=k$ one has
\begin{displaymath} 
\int_{S^{m-1}} \overline{M_t(\underline{\omega})} \ M_k(\underline{\omega}) \ w_k(\underline{\omega}) \ dS(\underline{\omega}) = 0 .
\end{displaymath}
(ii) Left outer Dunkl monogenics of different degree are orthogonal, i.e. for $Q_t \in {\mathcal M}^-_{\ell}(t)$ and $Q_k \in {\mathcal M}^-_{\ell}(k)$ with $t \not=k$ one has
\begin{displaymath} 
\int_{S^{m-1}} \overline{Q_t(\underline{\omega})} \ Q_k(\underline{\omega}) \ w_k(\underline{\omega}) \ dS(\underline{\omega}) = 0  .
\end{displaymath}
(iii) Any left inner and left outer Dunkl monogenic are orthogonal, i.e. for all $Q_t \in {\mathcal M}^-_{\ell}(t)$ and $M_k \in {\mathcal M}^+_{\ell}(k)$ one has
\begin{displaymath} 
\int_{S^{m-1}} \overline{Q_t(\underline{\omega})} \ M_k(\underline{\omega}) \ w_k(\underline{\omega}) \ dS(\underline{\omega}) = 0  .
\end{displaymath}
\end{theorem}
\textit{Proof.}\\
(i)-(ii) This follows immediately from the fact that Dunkl harmonics of different degree are orthogonal, i.e. \begin{displaymath}
\int_{S^{m-1}} H_k(\underline{\omega}) \ H_{\ell}(\underline{\omega}) \ w_k(\underline{\omega}) \ dS(\underline{\omega}) = 0
\end{displaymath}
if $k \not= \ell$ (see \cite{Dunkl}, p. 177).\\
(iii) The proof is based on the Cauchy theorem with $\Omega$ the closed unit ball $B(1)$ and hence $\partial \Omega$ the unit sphere $S^{m-1}$.\\
Take $Q_t \in {\mathcal M}^-_{\ell}(t)$ and $M_k \in {\mathcal M}^+_{\ell}(k)$. As at each point $\underline{\omega} \in S^{m-1}$, $n(\underline{\omega}) = \underline{\omega}$, we have that $d\sigma (\underline{\omega}) = \underline{\omega} \ dS(\underline{\omega})$ or $dS(\underline{\omega}) = -\underline{\omega} \ d\sigma (\underline{\omega})$. Hence we obtain
\begin{equation}\label{integraal}
\int_{S^{m-1}} \overline{Q_t(\underline{\omega})} \ M_k(\underline{\omega}) \ w_k(\underline{\omega}) \ dS(\underline{\omega})  =   \int_{S^{m-1}} \overline{\underline{\omega} \ Q_t(\underline{\omega})} \ d\sigma (\underline{\omega}) \ M_k(\underline{\omega}) \ w_k(\underline{\omega})  .
\end{equation}
As $Q_t \in {\mathcal M}^-_{\ell}(t)$, there exists $M_t \in {\mathcal M}^+_{\ell}(t)$ such that (see Proposition \ref{connectionmonogenic}): $M_t(\underline{\omega}) = \underline{\omega} \ Q_t(\underline{\omega})$. Hence equation (\ref{integraal}) becomes
\begin{displaymath}
\int_{S^{m-1}} \overline{Q_t(\underline{\omega})} \ M_k(\underline{\omega}) \ w_k(\underline{\omega}) \ dS(\underline{\omega})  = \int_{S^{m-1}} \overline{M_t(\underline{\omega})} \ d\sigma (\underline{\omega}) \ M_k(\underline{\omega}) \ w_k(\underline{\omega}) .
\end{displaymath}
Moreover, as $\overline{M_t}$ is right monogenic in $B(1)$, while $M_k$ is left monogenic in $B(1)$, the Cauchy theorem yields
\begin{displaymath}
\int_{S^{m-1}} \overline{Q_t(\underline{\omega})} \ M_k(\underline{\omega}) \ w_k(\underline{\omega}) \ dS(\underline{\omega})  =  0  . \ \ \ \square
\end{displaymath}
\section{Clifford-Gegenbauer polynomials on $B(1)$ related to the Dunkl Dirac operator}
\subsection{The Clifford-Gegenbauer polynomials on $B(1)$}\label{defCnops}
The Clifford-Gegenbauer polynomials on the unit ball $B(1)$ were introduced by Cnops (see \cite{Cnops}). They are orthogonal on $B(1)$ w.r.t. the weight function $(1-|\ux|^2)^{\alpha}$, $\alpha > -1$.
\begin{definition}~\\
Let $P_k$ be a spherical monogenic of degree $k$, $\alpha \in \mathbb{R}$, $\alpha > -1$ and $t$ a positive integer. Then
\begin{displaymath}
{\mathcal C}_t^{m,\alpha}(P_k) = D_{\alpha} D_{\alpha+1} \ldots D_{\alpha+t-1} \lbrack P_k \rbrack
\end{displaymath}
with
\begin{displaymath}
D_{\alpha} = (1-|\ux|^2) \partial_{\ux} - 2 (\alpha+1) \ux
\end{displaymath}
is a Clifford-Gegenbauer polynomial. Here $\partial_{\ux} = \sum_{j=1}^m e_j \partial_{x_j}$ is the so-called Dirac operator.
\end{definition}
\subsection{Definition Clifford-Gegenbauer polynomials on $B(1)$ related to the Dunkl Dirac}\label{defB1}
By analogy with the previous subsection, we first introduce for $\alpha \in \mathbb{R}$ and $\alpha>-1$ the operator
\begin{displaymath}
D_{\alpha}  =  (1-|\ux|^2) D_k - 2(\alpha +1) \ux  .
\end{displaymath}
The operator $D_{\alpha}$ can be rewritten as 
\begin{equation}\label{form2}
D_{\alpha} = (1-|\ux|^2)^{-\alpha} \ D_k (1-|\ux|^2)^{\alpha+1}  .
\end{equation}
Indeed, by means of the product rule (\ref{prod}), we obtain consecutively
\begin{eqnarray*}
& & (1-|\ux|^2)^{- \alpha} \ D_k \lbrack (1-|\ux|^2)^{\alpha + 1} \ f(\ux) \rbrack\\
& = & \sum_{i=1}^m \sum_A e_i (1-|\ux|^2)^{- \alpha} \ T_i \lbrack (1-|\ux|^2)^{\alpha+1} \ f_A(\ux) \rbrack \ e_A\\
& = & \sum_{i=1}^m \sum_A e_i (1-|\ux|^2)^{-\alpha} \ \lbrace \ T_i \lbrack (1-|\ux|^2)^{\alpha+1} \rbrack \ f_A(\ux) + (1-|\ux|^2)^{\alpha + 1} \ T_i \lbrack f_A(\ux) \rbrack \ \rbrace \ e_A\\
& = & \sum_{i=1}^m \sum_A e_i (1-|\ux|^2)^{-\alpha} \ \lbrace \ (\alpha + 1) (1-|\ux|^2)^{\alpha} \ (-2 x_i) \ f_A(\ux) + (1-|\ux|^2)^{\alpha + 1} \ T_i \lbrack f_A(\ux) \rbrack \ \rbrace \ e_A\\
& = & -2 (\alpha + 1) \ \ux \ f(\ux) + (1-|\ux|^2) \ D_k \lbrack f(\ux) \rbrack  .
\end{eqnarray*}

We now define the Clifford-Gegenbauer polynomials on $B(1)$ related to the Dunkl Dirac operator as follows:
\begin{definition}~\\
Let $M_k \in M_{\ell}^+(k)$ and $t$ a positive integer. Then
\begin{displaymath}
{\mathcal C}_{t,\mu}^{\alpha}(M_k)(\ux)  =  D_{\alpha} D_{\alpha+1} D_{\alpha+2} \ldots D_{\alpha+t-1} \lbrack M_k \rbrack
\end{displaymath}
is a Clifford-Gegenbauer polynomial on $B(1)$ of degree $t$ associated with $M_k$.
\end{definition}

Using Lemma \ref{belangrijk} we see that the precise form of the polynomials ${\mathcal C}_{t,\mu}^{\alpha}(M_k)$ depends only on the degree of the Dunkl-monogenic $M_k$, so we can write ${\mathcal C}_{t,\mu}^{\alpha}(M_k)(\ux) = {\mathcal C}_{t,\mu,k}^{\alpha}(\ux) \ M_k$.\\
The lower-degree polynomials take the following form:
\begin{eqnarray*}
{\mathcal C}_{0,\mu}^{\alpha}(M_k) & = & M_k\\
{\mathcal C}_{1,\mu}^{\alpha}(M_k) & = & -2(\alpha +1) \ux M_k\\
{\mathcal C}_{2,\mu}^{\alpha}(M_k) & = & 2 (\alpha +2) (2\alpha+2+2k+\mu) \ux^2 M_k + 2(\alpha +2) (2k+\mu) M_k.
\end{eqnarray*}
It is clear that ${\mathcal C}_{2s,\mu,k}^{\alpha}(\ux)$ will only contain even powers of $\ux$, while ${\mathcal C}_{2s+1,\mu,k}^{\alpha}(\ux)$ will only contain odd powers of $\ux$.
\subsection{Properties}\label{propertiesb1}
Using the definition, we immediately obtain that the Clifford-Gegenbauer polynomials on $B(1)$ satisfy the following recursion relation
\begin{eqnarray}\label{rec}
{\mathcal C}_{t+1,\mu}^{\alpha}(M_k) & = & D_{\alpha} \lbrack {\mathcal C}_{t,\mu}^{\alpha+1}(M_k) \rbrack\nonumber\\
& = & - 2 (\alpha +1) \ux \ {\mathcal C}_{t,\mu}^{\alpha+1}(M_k) + (1-|\ux|^2) \ D_k \lbrack {\mathcal C}_{t,\mu}^{\alpha+1}(M_k) \rbrack  .
\end{eqnarray}
The above result can be refined.
\begin{proposition}~\\
The Clifford-Gegenbauer polynomials satisfy the following recursion relations:
\begin{equation}\label{rec1}
{\mathcal C}_{2 \ell+1,\mu,k}^{\alpha}(\ux) = -2(\alpha+1) \ \ux \ {\mathcal C}_{2 \ell,\mu,k}^{\alpha+1}(\ux) + (1-|\ux|^2) \ D_k \lbrack {\mathcal C}_{2 \ell,\mu,k}^{\alpha+1}(\ux) \rbrack \ \ 
\end{equation}
and 
\begin{equation}\label{rec2}
{\mathcal C}_{2 \ell+2,\mu,k}^{\alpha}(\ux) = -2(\alpha+1) \ \ux \ {\mathcal C}_{2 \ell+1,\mu,k}^{\alpha+1}(\ux) + (1-|\ux|^2) \ \left( 2k \ \frac{\ux}{|\ux|^2} \ {\mathcal C}_{2 \ell+1,\mu,k}^{\alpha+1}(\ux) + D_k \lbrack {\mathcal C}_{2 \ell+1,\mu,k}^{\alpha+1}(\ux) \rbrack \right)  .
\end{equation}
\end{proposition}
\textit{Proof.} From (\ref{rec}) we immediately obtain (\ref{rec1}), whereas
\begin{displaymath}
{\mathcal C}^{\alpha}_{2\ell+2,\mu,k}(\ux) \ M_ k  =  -2 (\alpha+1) \ \ux \ {\mathcal C}^{\alpha+1}_{2\ell+1,\mu,k}(\ux) \ M_k + (1-|\ux|^2) \ D_k \lbrack {\mathcal C}^{\alpha+1}_{2\ell+1,\mu,k}(\ux) \ M_k  \rbrack  .
\end{displaymath}
Next by means of the product rule (\ref{prod}) and Proposition \ref{connectionmonogenic}, we have that 
\begin{eqnarray*}
D_k \lbrack {\mathcal C}^{\alpha+1}_{2\ell+1,\mu,k}(\ux) \ M_k  \rbrack &= & - D_k \left\lbrack r^{2k+\mu-1} \underline{\omega} \ {\mathcal C}^{\alpha+1}_{2\ell+1,\mu,k}(\ux) \ \underline{\omega} \ \frac{M_k}{r^{2k+\mu-1}} \right\rbrack\\
& = & - D_k \lbrack r^{2k+\mu-1} \ \underline{\omega} \ {\mathcal C}^{\alpha+1}_{2\ell+1,\mu,k}(\ux) \rbrack \ \underline{\omega} \ \frac{M_k}{r^{2k+\mu-1}} .
\end{eqnarray*}
Moreover, taking into account Proposition \ref{exprDirac} and the fact that the angular Dunkl Dirac operator only acts on the angular co-ordinates, yields
\begin{displaymath}
D_k \lbrack r^{2k+\mu-1} \ \underline{\omega} \ {\mathcal C}^{\alpha+1}_{2\ell+1,\mu,k}(\ux) \rbrack 
 =  - (2k+\mu-1) \ r^{2k+\mu-2} \ {\mathcal C}^{\alpha+1}_{2\ell+1,\mu,k}(\ux) + \underline{\omega} \ r^{2k+\mu-1} \underline{\omega} \ \partial_r \lbrack {\mathcal C}^{\alpha+1}_{2\ell+1,\mu,k}(\ux) \rbrack  .
\end{displaymath}
Furthermore, as
\begin{displaymath}
\underline{\omega} \ \partial_r \lbrack {\mathcal C}^{\alpha+1}_{2\ell+1,\mu,k}(\ux) \rbrack = \left( D_k - \frac{\underline{\omega}}{r} \Gamma_k \right) \lbrack {\mathcal C}^{\alpha+1}_{2\ell+1,\mu,k}(\ux) \rbrack = D_k \left\lbrack {\mathcal C}^{\alpha+1}_{2\ell+1,\mu,k}(\ux) \right\rbrack - \frac{(\mu -1) \underline{\omega}}{r} \ {\mathcal C}^{\alpha+1}_{2\ell+1,\mu,k}(\ux) ,
\end{displaymath}
we finally obtain
\begin{eqnarray*}
& & {\mathcal C}^{\alpha}_{2\ell+2,\mu,k}(\ux)\\
& = & -2 (\alpha +1 ) \ux \ {\mathcal C}^{\alpha+1}_{2\ell+1,\mu,k}(\ux) - (1-|\ux|^2) \biggl\lbrace - (2k+\mu-1)\  \frac{{\mathcal C}^{\alpha+1}_{2\ell+1,\mu,k}(\ux)}{r} \ \underline{\omega}\\
& & + \underline{\omega} \ D_k \lbrack {\mathcal C}^{\alpha+1}_{2\ell+1,\mu,k}(\ux) \rbrack \ \underline{\omega} + \frac{(\mu-1)}{r} \ \underline{\omega} \ {\mathcal C}^{\alpha+1}_{2\ell+1,\mu,k}(\ux) \biggr\rbrace\\
& = & -2 (\alpha +1 ) \ux \ {\mathcal C}^{\alpha+1}_{2\ell+1,\mu,k}(\ux) + (1-|\ux|^2) \ \left( 2k \frac{\ux}{|\ux|^2} \ {\mathcal C}^{\alpha+1}_{2\ell+1,\mu,k}(\ux) + D_k \lbrack {\mathcal C}^{\alpha+1}_{2\ell+1,\mu,k}(\ux) \rbrack \right) . \ \ \square
\end{eqnarray*}

Taking into account (\ref{form2}), it is easily seen that there also exists a Rodrigues formula.
\begin{theorem}[Rodrigues formula]~\\
The Clifford-Gegenbauer polynomials on $B(1)$ take the form
\begin{displaymath}
{\mathcal C}_{t,\mu}^{\alpha}(M_k) = (1-|\ux|^2)^{-\alpha} \ D_k^{t}  \biggl\lbrack (1-|\ux|^2)^{\alpha+t} \ M_k \biggr\rbrack  .
\end{displaymath}
\end{theorem}
Moreover, the Clifford-Gegenbauer polynomials on $B(1)$ satisfy an annihilation equation.
\begin{theorem}[Annihilation equation]~\\ ${\mathcal C}_{t,\mu}^{\alpha}(M_k)$ satisfies
\begin{displaymath}
D_k \lbrack {\mathcal C}_{t,\mu}^{\alpha}(M_k) \rbrack = C(\alpha,t,\mu,k) \ {\mathcal C}_{t-1,\mu}^{\alpha+1}(M_k)
\end{displaymath}
with
\begin{equation*}
C(\alpha,t,\mu,k) = 
  \begin{cases}
       t(2 \alpha+t+\mu+2k) & \text{if $t$ is even},\\
       (2 \alpha+t+1)(t+\mu+2k-1) & \text{if $t$ is odd}.
  \end{cases}
\end{equation*}
\end{theorem}
\textit{Proof.} Let us write down the expansion of the Clifford-Gegenbauer polynomials:
\begin{displaymath}
{\mathcal C}_{2t,\mu}^{\alpha}(M_k)(\ux) = \sum_{i=0}^t a_{2i}^{2t,\alpha} \ \ux^{2i} \ M_k \qquad \mathrm{and} \qquad {\mathcal C}_{2t+1,\mu}^{\alpha}(M_k)(\ux) = \sum_{i=0}^t a_{2i+1}^{2t+1,\alpha} \ \ux^{2i+1} \ M_k  .
\end{displaymath}
From the recursion formula (\ref{rec}) and Lemma \ref{belangrijk}, we know that the following relations between the coefficients hold:
\begin{equation}\label{geld1}
a_{2i}^{2t,\alpha} = -(2i+2k+\mu) \ a_{2i+1}^{2t-1,\alpha+1} - (2i+2k+2\alpha+\mu) \ a_{2i-1}^{2t-1,\alpha+1}
\end{equation}
and
\begin{equation}\label{geld2}
a_{2i+1}^{2t+1,\alpha} = -2(\alpha+1+i) \ a_{2i}^{2t,\alpha+1} - (2i+2) \ a_{2i+2}^{2t,\alpha+1}  .
\end{equation}
We need to prove
\begin{equation}\label{tebew1}
-2i \ a_{2i}^{2t,\alpha} = 2t \ (2 \alpha+2t+\mu + 2k) \ a_{2i-1}^{2t-1,\alpha+1}
\end{equation}
and
\begin{equation}\label{tebew2}
-(2i+2k+\mu) \ a_{2i+1}^{2t+1,\alpha} = (2 \alpha + 2t +2) (2t+\mu+2k) \ a_{2i}^{2t,\alpha+1}  .
\end{equation}
It is easy to check that the theorem holds for $t=0,1$. Using (\ref{geld1}) and (\ref{geld2}) the theorem can then be proved by induction on $t$. \ \ \ \ \ $\square$

By acting on the above annihilation equation with $D_{\alpha}$ we immediately obtain the differential equation satisfied by the Clifford-Gegenbauer polynomials.
\begin{theorem}[Differential equation]~\\
${\mathcal C}_{t,\mu}^{\alpha}(M_k)$ is a solution of the following differential equation:
\begin{displaymath}
 (1-|\ux|^2) \Delta_k \lbrack {\mathcal C}_{t,\mu}^{\alpha}(M_k) \rbrack + 2 (\alpha+1) \ux D_k \lbrack {\mathcal C}_{t,\mu}^{\alpha}(M_k) \rbrack + C(\alpha,t,\mu,k) {\mathcal C}_{t,\mu}^{\alpha}(M_k) = 0 .
\end{displaymath}
\end{theorem}
The above equation should be compared with the classical differential equation of the Gegenbauer polynomials on the real line
\begin{displaymath}
(1-x^2) \frac{d^2}{dx^2} C_n^{\lambda}(x) - (2 \lambda +1) \ x \frac{d}{dx} C_n^{\lambda}(x) + n (n+2 \lambda) \ C_n^{\lambda}(x) = 0
\end{displaymath}
where it should be noticed that $\lambda = \alpha +\frac{1}{2}$.\\
\\
Moreover, combining the annihilation equation and the recursion formula (\ref{rec}) we obtain the following recurrence relation.
\begin{theorem}[Recurrence relation]~\\
${\mathcal C}_{t,\mu}^{\alpha}(M_k)$ satisfies the recurrence relation:
\begin{displaymath}
{\mathcal C}_{t+1,\mu}^{\alpha}(M_k)(\ux) + 2(\alpha + 1) \ \ux \ {\mathcal C}_{t,\mu}^{\alpha+1}(M_k) - C(\alpha+1,t,\mu,k) \ (1-|\ux|^2) \ {\mathcal C}_{t-1,\mu}^{\alpha+2}(M_k) = 0  .
\end{displaymath}
\end{theorem}
Note that this recursion formula is the Dunkl analogon of the classical one-dimensional Gegenbauer recurrence relation:
\begin{displaymath}
(n+1) \ C_{n+1}^{\lambda}(x) - (n+2\lambda) \ x \ C_n^{\lambda}(x) + 2 \lambda (1-x^2) \ C_{n-1}^{\lambda+1}(x) = 0  .
\end{displaymath}

It is now possible to express the Clifford-Gegenbauer polynomials on $B(1)$ in terms of the Jacobi polynomials on the real line.
\begin{theorem}[Closed form]\label{closed}~\\
The Clifford-Gegenbauer polynomials on $B(1)$ can be written in terms of the Jacobi polynomials on the real line as
\begin{eqnarray}\label{explicit1}
{\mathcal C}_{2t,\mu,k}^{\alpha}(\ux) & = & 2^{2t} \ (\alpha+t+1)_t \ t! \ P_t^{(\mu/2+k-1,\alpha)}(1+2 \ux^2)\\
& = & (-1)^t \ 2^{2t} \ (\alpha+t+1)_t \ t! \ P_t^{(\alpha,\mu/2+k-1)}(2 |\ux|^2-1)\nonumber
\end{eqnarray}
and
\begin{eqnarray}\label{explicit2}
{\mathcal C}_{2t+1,\mu,k}^{\alpha}(\ux) & = & -2^{2t+1} \ (\alpha+t+1)_{t+1} \ t! \ \ux \ P_t^{(\mu/2+k,\alpha)}(1+2 \ux^2)\\
& = & (-1)^{t+1} \ 2^{2t+1} \ (\alpha+t+1)_{t+1} \ t! \ \ux \ P_t^{(\alpha, \mu/2+k)}(2 |\ux|^2 -1)\nonumber
\end{eqnarray}
with $(a)_p = a (a+1) \ldots (a+p-1) = \displaystyle{\frac{\Gamma(a+p)}{\Gamma(a)}}$ the Pochhammer symbol and
\begin{displaymath}
P_t^{(\alpha,\beta)}(x) = \frac{\Gamma(\alpha+t+1)}{t! \ \Gamma(\alpha+\beta+t+1)} \ \sum_{i=0}^t {t \choose i} \frac{\Gamma(\alpha+\beta+t+i+1)}{\Gamma(\alpha+i+1)} \ \left(\frac{x-1}{2} \right)^i  .
\end{displaymath}
\end{theorem}
\textit{Proof.} Recall that
\begin{displaymath}
{\mathcal C}_{2t,\mu}^{\alpha}(M_k)(\ux) = \sum_{i=0}^t a_{2i}^{2t,\alpha} \ \ux^{2i} \ M_k \ \ .
\end{displaymath}
By means of the annihilation equations (\ref{tebew1}) and (\ref{tebew2}) in terms of the coefficients $a_i^{s,\alpha}$, we can write $a_{2i}^{2t,\alpha}$ in terms of $a_0^{2t-2i,\alpha+2i}$ :
\begin{eqnarray}\label{res1}
a_{2i}^{2t,\alpha} & = & - \frac{t}{i} \ (2 \alpha+2t+\mu+2k) \ a_{2i-1}^{2t-1,\alpha+1}\nonumber\\
& = & (-1)^2 2^2 \ \frac{t}{i} \ \left(\alpha+t+\frac{\mu}{2} +k \right) \ (\alpha+t+1) \ \frac{\left( t + \frac{\mu}{2} +k-1 \right)}{\left( i+k+\frac{\mu}{2}-1 \right)} \ a_{2i-2}^{2t-2,\alpha+2}\nonumber\\
& = & ...\nonumber\\
& = & 2^{2i} \ { t \choose i} \ \left(\alpha+t+\frac{\mu}{2}+k \right)_i \ (\alpha+t+1)_i \ \frac{\Gamma \left( t+ \frac{\mu}{2} +k \right)}{\Gamma \left( t+\frac{\mu}{2} +k-i \right)}\nonumber\\
& & \times \frac{\Gamma \left( k+ \frac{\mu}{2} \right)}{\Gamma \left( k+ \frac{\mu}{2} +i \right)} \ a_0^{2t-2i,\alpha+2i}  .
\end{eqnarray}
Now we look for an expression of $a_0^{2t,\alpha}$. By means of successively (\ref{geld1}) and (\ref{tebew2}) we are able to write $a_0^{2t,\alpha}$ in terms of $a_0^{0,\alpha + 2t} = 1$ :
\begin{eqnarray}\label{res2}
a_0^{2t,\alpha} & = & -(2k+\mu) a_1^{2t-1,\alpha+1}\nonumber\\
& = & 2^2 \ (\alpha+1+t)(t-1+\mu/2 +k) \ a_0^{2t-2,\alpha+2}\nonumber\\
& = & 2^4 \ (\alpha+1+t)(\alpha+t+2)(t+\mu/2+k-1)(t+\mu/2+k-2) \ a_0^{2t-4,\alpha+4}\nonumber\\
& = & ...\nonumber\\
& = & 2^{2t} \ (\alpha+t+1)_t \ \frac{\Gamma(\mu/2+k+t)}{\Gamma (\mu/2+k)} \ a_0^{0,\alpha+2t}  .
\end{eqnarray}
Combining (\ref{res1}) and (\ref{res2}) yields
\begin{displaymath}
a_{2i}^{2t,\alpha}  = 2^{2t} \ {t \choose i} \frac{\Gamma(t+\mu/2+k)}{\Gamma(k+\mu/2+i)} \ (\alpha+t+\mu/2+k)_i \ (\alpha+t+1)_t  ,
\end{displaymath}
from which we indeed obtain (\ref{explicit1}).\\
The formula for $a_{2i+1}^{2t+1,\alpha}$ now follows from the annihilation equation (\ref{tebew2}):
\begin{eqnarray*}
a_{2i+1}^{2t+1,\alpha} & =& - 2(\alpha+t+1) \frac{t+\mu/2+k}{i+k+\mu/2} \ a_{2i}^{2t,\alpha+1}\\
& = & - 2^{2t+1} {t \choose i} \ \frac{\Gamma(t+\mu/2+k+1)}{\Gamma(k+\mu/2+i+1)} \ (\alpha+t+\mu/2+k+1)_i \ (\alpha+t+1)_{t+1} ,
\end{eqnarray*}
which leads to (\ref{explicit2}). \ \ \ \ $\square$\vspace{0,2cm}\\
Let us mention the following corollary of the previous theorem.
\begin{corollary}~\\
The Clifford-Gegenbauer polynomials on $B(1)$ satisfy
\begin{displaymath}
{\mathcal C}^{\alpha}_{2t+1,\mu,k}(\ux) = -2 (\alpha+2t+1) \ \ux \ {\mathcal C}^{\alpha}_{2t,\mu,k+1}(\ux) .
\end{displaymath}
\end{corollary}
By means of Theorem \ref{closed} we are able to prove the Dunkl-analogon of the classical one-dimensional three-term Gegenbauer recurrence relation:
\begin{displaymath}
n \ C_n^{\lambda}(x) = 2(n+\lambda-1) \ x \ C_{n-1}^{\lambda}(x) - (n+2\lambda-2) \ C_{n-2}^{\lambda}(x) \qquad , \qquad n=2,3,\ldots
\end{displaymath}
\begin{theorem}[Three-term recurrence relation]~\\
$C_{t,\mu,k}^{\alpha}(\ux)$ satisfies the three-term recurrence relation:
\begin{displaymath}
\frac{D(\alpha,t,\mu,k)}{2(\alpha+t)} \ C_{t,\mu,k}^{\alpha}(\ux) = - \left( \alpha+ \frac{\mu}{2} +k+t-1 \right) \ \ux \ C_{t-1,\mu,k}^{\alpha}(\ux) + (\alpha+t-1) \ E(t,\mu,k) \ C_{t-2,\mu,k}^{\alpha}(\ux)
\end{displaymath}
with
\begin{equation*}
D(\alpha,t,\mu,k) = 
  \begin{cases}
       \alpha + \frac{t}{2} & \text{if $t$ is even},\\
       \alpha+\frac{\mu}{2}+k+\frac{t}{2}-\frac{1}{2} & \text{if $t$ is odd}
  \end{cases}
\end{equation*}
and
\begin{equation*}
E(t,\mu,k) = 
  \begin{cases}
       \mu+2k-2+t & \text{if $t$ is even},\\
       t-1 & \text{if $t$ is odd}.
  \end{cases}
\end{equation*}
\end{theorem}
\textit{Proof.} This result follows from the following contiguous relations of the classical Jacobi polynomials on the real line (see \cite{Magnus}):
\begin{displaymath}
\left( \frac{\alpha}{2} + \frac{\beta}{2} + \ell +1 \right) \ (1+y) \ P_{\ell}^{(\alpha,\beta+1)}(y) = (\beta + \ell+1) \ P_{\ell}^{(\alpha,\beta)}(y) + (\ell+1) \ P_{\ell+1}^{(\alpha,\beta)}(y)
\end{displaymath}
and 
\begin{displaymath}
(\alpha +\beta+2n) \ P_n^{(\alpha,\beta-1)}(y) = (\alpha+\beta+n) \ P_n^{(\alpha,\beta)}(y) + (\alpha+n) \ P_{n-1}^{(\alpha,\beta)}(y)  . \ \ \ \ \square
\end{displaymath}
\subsection{Orthogonality}
Let us consider the inner product:
\begin{displaymath}
\langle f,g \rangle_{\alpha} = \int_{B(1)} \overline{f(\ux)} \ g(\ux) \ (1-|\ux|^2)^{\alpha} \ w_k(\ux) \ dV(\ux).
\end{displaymath}
\begin{proposition}\label{prop2}~\\
The operators $D_{\alpha}$ and $D_k$ are dual with respect to $\langle . \ , \ .\rangle_{\alpha}$, i.e. for all $f,g \in C_1 \bigl( B(1) \bigr)$
\begin{displaymath}
\langle D_{\alpha} \lbrack f \rbrack \ , \ g \rangle_{\alpha} \ = \  \langle f \ , \ D_k \lbrack g \rbrack \rangle_{\alpha +1}  .
\end{displaymath}
\end{proposition}
\textit{Proof.} By means of Theorem \ref{stokes}, we find
\begin{eqnarray*}
 \langle D_{\alpha} \lbrack f \rbrack \ , \ g \rangle_{\alpha} & = & - \int_{B(1)} \lbrack (1-|\ux|^2)^{\alpha+1} \ \overline{f} \rbrack D_k \ g(\ux) \ w_k(\ux) \ dV(\ux)\\
& = & - \biggl( \int_{\partial B(1)} (1-|\ux|^2)^{\alpha+1} \ \overline{f} \ w_k(\ux) \ d\sigma(\ux) \ g(\ux)\\
& & - \int_{B(1)} (1-|\ux|^2)^{\alpha+1} \ \overline{f} \ D_k \lbrack g(\ux) \rbrack \ w_k(\ux) \ dV(\ux) \biggr)\\
& = & \int_{B(1)} \overline{f} \ D_k \lbrack g \rbrack \ (1-|\ux|^2)^{\alpha+1} \ w_k(\ux) \ dV(\ux) \ = \ \langle f \ , \ D_k \lbrack g \rbrack \rangle_{\alpha+1}  ,
\end{eqnarray*}
since the surface term clearly vanishes. $\square$\\

Using the above proposition, we are now able to prove an orthogonality relation for the Clifford-Gegenbauer polynomials on $B(1)$.
\begin{theorem}[Orthogonality relation]~\\
If $ s \not= t$ or $k \not= \ell$, then
\begin{displaymath}
\langle {\mathcal C}_{t,\mu}^{\alpha}(M_k), {\mathcal C}_{s,\mu}^{\alpha}(M_{\ell})\rangle_{\alpha} = 0  .
\end{displaymath}
\end{theorem}
\textit{Proof.}\\
a) Suppose that $s \not= t$ and $t>s$ (the case where $t<s$ is similar). We obtain
\begin{eqnarray*}
 \langle {\mathcal C}_{t,\mu}^{\alpha}(M_k), {\mathcal C}_{s,\mu}^{\alpha}(M_{\ell})\rangle_{\alpha} & = & \langle D_{\alpha} \ D_{\alpha+1} \ldots D_{\alpha +t-1} \lbrack M_k \rbrack , {\mathcal C}_{s,\mu}^{\alpha}(M_{\ell})\rangle_{\alpha}\\
& = & \langle D_{\alpha+1} \ldots D_{\alpha +t-1} \lbrack M_k \rbrack , D_k \lbrack {\mathcal C}_{s,\mu}^{\alpha}( M_{\ell}) \rbrack \rangle_{\alpha+1}\\
& = & . \ . \ . \ .\\
& = & \langle M_k , D_k^t \lbrack {\mathcal C}_{s,\mu}^{\alpha} (M_{\ell}) \rbrack \rangle_{\alpha+t} \ = \ 0  .
\end{eqnarray*}
b) If $s=t$ and $k \not= \ell$, the result follows from the fact that the left inner Dunkl monogenics of different degree are orthogonal (see Theorem \ref{orthomon}). \ \ \ \ $\square$\\

Using the orthonormality relation of the classical Jacobi polynomials, we are able to calculate the normalization constants.
\begin{lemma}~\\
We have that
\begin{multline*}
 \left< {\mathcal C}_{2t,\mu}^{\alpha}(M_k) \ , \ {\mathcal C}_{2t,\mu}^{\alpha}(M_k) \right>_{\alpha}\\
 =  2^{4t-1} \ t! \ \frac{(\alpha +t +1)_t \ \Gamma(\alpha + 2t +1)}{\left( \frac{\mu}{2} + k +t \right)_{\alpha} \ \left( \frac{\mu}{2} + k+\alpha+2t \right)} \ \int_{S^{m-1}} \overline{M_k(\underline{\omega})} \ M_k(\underline{\omega}) \ w_k(\underline{\omega}) \ dS(\underline{\omega})
\end{multline*}
\begin{multline*}
 \left< {\mathcal C}_{2t+1,\mu}^{\alpha}(M_k) \ , \ {\mathcal C}_{2t+1,\mu}^{\alpha}(M_k)\right>_{\alpha}\\
 =  - 2^{4t+1} \ t! \ \frac{(\alpha +t +1)_{t+1} \ \Gamma(\alpha + 2t +2)}{\left( \frac{\mu}{2} + k +t+1 \right)_{\alpha} \ \left( \frac{\mu}{2} + k+\alpha+2t+1 \right)} \ \int_{S^{m-1}} \overline{M_k(\underline{\omega})} \ M_k(\underline{\omega}) \ w_k(\underline{\omega}) \ dS(\underline{\omega})  .
\end{multline*}
\end{lemma}\vspace{0,5cm}
\begin{remark}[Gegenbauer polynomials on $B(1)$ related to the Dunkl Laplacian]
Note that if we calculate $D_{\alpha} D_{\alpha +1}$ we obtain
\begin{eqnarray*}
{\mathcal D}_{\alpha} & := & D_{\alpha} \ D_{\alpha +1}\\
& = & - (1-|\ux|^2)^2 \ \Delta_k - 2(\alpha+2)(2 \alpha+2+\mu) \ |\ux|^2 + 4(\alpha+2) \ (1-|\ux|^2) \ \mathbb{E} + 2 (\alpha +2) \mu  ,
\end{eqnarray*}
where we have used the product rule (\ref{prod}) and the relation (\ref{fundamenteel}).\\
As this operator is scalar, it makes sense to let it act on a Dunkl harmonic $H_k$ instead of a Dunkl monogenic. Hence, we can define
\begin{displaymath}
{\mathcal C}_{2t}^{\mu,\alpha}(H_k) = {\mathcal D}_{\alpha} {\mathcal D}_{\alpha+2}{\mathcal D}_{\alpha+4} \ldots {\mathcal D}_{\alpha+2t-2} \lbrack H_k \rbrack  .
\end{displaymath}
For these Gegenbauer polynomials associated with $H_k$ we can also derive a Rodrigues formula, a differential equation, an annihilation equation and a recurrence relation. Moreover, in terms of the Jacobi polynomials on the real line, they take the following form
\begin{displaymath}
{\mathcal C}_{2t}^{\mu,\alpha}(H_k) = 2^{2t} \ (\alpha + t+1)_t \ t! \ P_t^{\left( \frac{\mu}{2}+k-1,\alpha \right)} (1+2 \ux^2) \ H_k  .
\end{displaymath}
Note that special cases of these scalar polynomials and there associated weights have already been studied in \cite{yuan1, yuan2}.\end{remark}
\section{Clifford-Gegenbauer polynomials on $\mathbb{R}^m$ related to the Dunkl Dirac operator}
\subsection{The generalized Clifford-Gegenbauer polynomials}
The generalized Clifford-Gegenbauer polynomials (see e.g. \cite{WCAA7} and \cite{wave}) are defined as follows:
\begin{definition}~\\
Let $P_k$ be a spherical monogenic of degree $k$, $\alpha \in \mathbb{R}$ and $t$ a positive integer. Then
\begin{displaymath}
{\mathcal G}^{m,\alpha}_t(P_k) = D_{\alpha} D_{\alpha +1} \ldots D_{\alpha+t-1} \lbrack P_k \rbrack
\end{displaymath}
with 
\begin{displaymath}
D_{\alpha} = (1+|\ux|^2)\partial_{\ux} + 2 (\alpha+1) \ux
\end{displaymath}
is a Clifford-Gegenbauer polynomial.
\end{definition}
This second type of Clifford-Gegenbauer polynomials came into play while studying wavelets in the Clifford analysis setting. The polynomials, originally defined in a completely different way as by Cnops in \cite{Cnops}, are the desired building blocks for new higher dimensional wavelet kernels (see \cite{wave}). They satisfy certain orthogonality relations on the whole of $\mathbb{R}^m$ w.r.t. the weight function $(1+|\ux|^2)^{\alpha}$, $\alpha \in \mathbb{R}$.
\subsection{Definition Clifford-Gegenbauer polynomials on $\mathbb{R}^m$ related to the Dunkl Dirac}
By analogy with the previous subsection, we first introduce for $\alpha \in \mathbb{R}$ the operator
\begin{displaymath}
D_{\alpha}  =  (1+|\ux|^2) D_k + 2(\alpha +1) \ux  .
\end{displaymath}
Similarly as in subsection \ref{defB1} one can verify that the operator $D_{\alpha}$ can be rewritten as 
\begin{equation}\label{alternform}
D_{\alpha} = (1+|\ux|^2)^{-\alpha} \ D_k (1+|\ux|^2)^{\alpha+1}  .
\end{equation}

We now define the Clifford-Gegenbauer polynomials on $\mathbb{R}^m$ related to the Dunkl Dirac operator as follows:
\begin{definition}~\\
Let $M_k \in M_{\ell}^+(k)$ and $t$ a positive integer. Then
\begin{displaymath}
{\mathcal G}_{t,\mu}^{\alpha}(M_k)(\ux)  =  D_{\alpha} D_{\alpha+1} D_{\alpha+2} \ldots D_{\alpha+t-1} \lbrack M_k \rbrack
\end{displaymath}
is a Clifford-Gegenbauer polynomial on $\mathbb{R}^m$ of degree $t$ associated with $M_k$.
\end{definition}

Again using Lemma \ref{belangrijk}, it is clear that the precise form of the polynomials ${\mathcal G}_{t,\mu}^{\alpha}(M_k)$ depends only on the degree of the Dunkl-monogenic $M_k$, so we can write ${\mathcal G}_{t,\mu}^{\alpha}(M_k)(\ux) = {\mathcal G}_{t,\mu,k}^{\alpha}(\ux) \ M_k$.\\
The lower-degree polynomials take the following form:
\begin{eqnarray*}
{\mathcal G}_{0,\mu}^{\alpha}(M_k) & = & M_k\\
{\mathcal G}_{1,\mu}^{\alpha}(M_k) & = & 2(\alpha +1) \ux M_k\\
{\mathcal G}_{2,\mu}^{\alpha}(M_k) & = & 2 (\alpha +2) (2\alpha+2+2k+\mu) \ux^2 M_k - 2(\alpha +2) (2k+\mu) M_k.
\end{eqnarray*}
It is clear that ${\mathcal G}_{2s,\mu,k}^{\alpha}(\ux)$ will only contain even powers of $\ux$, while ${\mathcal G}_{2s+1,\mu,k}^{\alpha}(\ux)$ will only contain odd powers of $\ux$.
\subsection{Properties}
In this subsection we collect the properties of the Clifford-Gegenbauer polynomials on $\mathbb{R}^m$. As the proofs are similar as those in subsection \ref{propertiesb1}, we omit them.

Using the definition, we immediately obtain that the Clifford-Gegenbauer polynomials on $\mathbb{R}^m$ satisfy the following recursion relation
\begin{eqnarray}\label{recrm}
{\mathcal G}_{t+1,\mu}^{\alpha}(M_k) & = & D_{\alpha} \lbrack {\mathcal G}_{t,\mu}^{\alpha+1}(M_k) \rbrack\nonumber\\
& = & 2 (\alpha +1) \ux \ {\mathcal G}_{t,\mu}^{\alpha+1}(M_k) + (1+|\ux|^2) \ D_k \lbrack {\mathcal G}_{t,\mu}^{\alpha+1}(M_k) \rbrack  .
\end{eqnarray}
In particular, we obtain the following
\begin{proposition}~\\
The Clifford-Gegenbauer polynomials satisfy the following recursion relations:
\begin{equation}
{\mathcal G}_{2 \ell+1,\mu,k}^{\alpha}(\ux) = 2(\alpha+1) \ \ux \ {\mathcal G}_{2 \ell,\mu,k}^{\alpha+1}(\ux) + (1+|\ux|^2) \ D_k \lbrack {\mathcal G}_{2 \ell,\mu,k}^{\alpha+1}(\ux) \rbrack \ \ 
\end{equation}
and 
\begin{equation}
{\mathcal G}_{2 \ell+2,\mu,k}^{\alpha}(\ux) = 2(\alpha+1) \ \ux \ {\mathcal G}_{2 \ell+1,\mu,k}^{\alpha+1}(\ux) + (1+|\ux|^2) \ \left( 2k \ \frac{\ux}{|\ux|^2} \ {\mathcal G}_{2 \ell+1,\mu,k}^{\alpha+1}(\ux) + D_k \lbrack {\mathcal G}_{2 \ell+1,\mu,k}^{\alpha+1}(\ux) \rbrack \right)  .
\end{equation}
\end{proposition}

Taking into account (\ref{alternform}), it is easily seen that there also exists a Rodrigues formula.
\begin{theorem}[Rodrigues formula]~\\
The Clifford-Gegenbauer polynomials on $\mathbb{R}^m$ take the form
\begin{displaymath}
{\mathcal G}_{t,\mu}^{\alpha}(M_k) = (1+|\ux|^2)^{-\alpha} \ D_k^{t}  \biggl\lbrack (1+|\ux|^2)^{\alpha+t} \ M_k \biggr\rbrack  .
\end{displaymath}
\end{theorem}
Moreover, the Clifford-Gegenbauer polynomials on $\mathbb{R}^m$ satisfy an annihilation equation.
\begin{theorem}[Annihilation equation]~\\ ${\mathcal G}_{t,\mu}^{\alpha}(M_k)$ satisfies
\begin{displaymath}
D_k \lbrack {\mathcal G}_{t,\mu}^{\alpha}(M_k) \rbrack = - C(\alpha,t,\mu,k) \ {\mathcal G}_{t-1,\mu}^{\alpha+1}(M_k)  .
\end{displaymath}
\end{theorem}
From the above annihilation equation we obtain the differential equation
\begin{theorem}[Differential equation]~\\
${\mathcal G}_{t,\mu}^{\alpha}(M_k)$ is a solution of the following differential equation:
\begin{displaymath}
 (1+|\ux|^2) \Delta_k \lbrack {\mathcal G}_{t,\mu}^{\alpha}(M_k) \rbrack - 2 (\alpha+1) \ux D_k \lbrack {\mathcal G}_{t,\mu}^{\alpha}(M_k) \rbrack - C(\alpha,t,\mu,k) {\mathcal G}_{t,\mu}^{\alpha}(M_k) = 0  .
\end{displaymath}
\end{theorem}

Moreover, combining the annihilation equation and the recursion formula (\ref{recrm}) we find the following recurrence relation.
\begin{theorem}[Recurrence relation]~\\
${\mathcal G}_{t,\mu}^{\alpha}(M_k)$ satisfies the recurrence relation:
\begin{displaymath}
{\mathcal G}_{t+1,\mu}^{\alpha}(M_k)(\ux) - 2(\alpha + 1) \ \ux \ {\mathcal G}_{t,\mu}^{\alpha+1}(M_k) + C(\alpha+1,t,\mu,k) \ (1+|\ux|^2) \ {\mathcal G}_{t-1,\mu}^{\alpha+2}(M_k) = 0  .
\end{displaymath}
\end{theorem}

We can now express the Clifford-Gegenbauer polynomials on $\mathbb{R}^m$ in terms of the Jacobi polynomials on the real line.
\begin{theorem}[Closed form]\label{closed2}~\\
The Clifford-Gegenbauer polynomials on $\mathbb{R}^m$ can be written in terms of the Jacobi polynomials on the real line as
\begin{eqnarray*}
{\mathcal G}_{2t,\mu,k}^{\alpha}(\ux) & = & (-1)^t \ 2^{2t} \ (\alpha+t+1)_t \ t! \ P_t^{(\mu/2+k-1,\alpha)}(1-2 \ux^2)\\
{\mathcal G}_{2t+1,\mu,k}^{\alpha}(\ux) & = & (-1)^t \ 2^{2t+1} \ (\alpha+t+1)_{t+1} \ t! \ \ux \ P_t^{(\mu/2+k,\alpha)}(1-2 \ux^2)  .
\end{eqnarray*}
\end{theorem}
Finally, the previous theorem enables us to prove the following result.
\begin{theorem}[Three-term recurrence relation]~\\
${\mathcal G}_{t,\mu,k}^{\alpha}(\ux)$ satisfies the three-term recurrence relation:
\begin{displaymath}
\frac{D(\alpha,t,\mu,k)}{2(\alpha+t)} \ {\mathcal G}_{t,\mu,k}^{\alpha}(\ux) = \left( \alpha + \frac{\mu}{2} +k+t-1 \right) \ \ux \ {\mathcal G}_{t-1,\mu,k}^{\alpha}(\ux) - (\alpha +t-1) \ E(t,\mu,k) \ {\mathcal G}_{t-2,\mu,k}^{\alpha}(\ux)  .
\end{displaymath}
\end{theorem}
\subsection{Orthogonality}
The orthogonality of the Clifford-Gegenbauer polynomials on $\mathbb{R}^m$ must be treated in a completely different way than the orthogonality of their counterparts on the unit ball $B(1)$. It must be expressed in terms of a bilinear form instead of an integral (as was also done in \cite{DBS3,WCAA7} for the super resp. classical case).

Let us start by computing the following integral on the Euclidean space $\mathbb{R}^m$ with $M_k$ and $M_{\ell}$ left solid inner Dunkl monogenics of order $k$, respectively $\ell$ :
\begin{eqnarray}\label{uitgewerkteint}
& & \hspace{-1,5 cm} \int_{\mathbb{R}^m} \overline{\ux^s \ M_k(\ux)} \ \ux^t \ M_{\ell}(\ux) \ (1+|\ux|^2)^{\alpha} \ w_k(\ux) \ dV(\ux)\nonumber\\
& = & \int_0^{+ \infty} r^{s+t+k+\ell+\mu-1} \ (1+r^2)^{\alpha} \ dr \ \int_{S^{m-1}} \overline{M_k(\underline{\omega})} \ \overline{\underline{\omega}^s} \ \underline{\omega}^t \ M_{\ell}(\underline{\omega}) \ w_k(\underline{\omega}) \ dS(\underline{\omega})\nonumber\\
&= & \frac{1}{2} \ B \left( \frac{s+t+k+\ell+\mu}{2} , - \left( \frac{s+t+k+\ell+\mu}{2} \right) - \alpha \right)\nonumber\\ 
& & \int_{S^{m-1}} \overline{M_k(\underline{\omega})} \ \overline{\underline{\omega}^s} \ \underline{\omega}^t \ M_{\ell}(\underline{\omega}) \ w_k(\underline{\omega}) \ dS(\underline{\omega})
\end{eqnarray}
using the Beta function
\begin{displaymath}
B(x,y) =  \int_0^{+ \infty} u^{x-1} \ (1+u)^{-x-y} \ du \ \ , \ \ \mathrm{Re}(x) >0 \ \ , \ \ \mathrm{Re}(y) > 0  .
\end{displaymath}
Naturally, the last equality in (\ref{uitgewerkteint}) only holds if $2 \alpha < -(k+s+t+\ell+\mu)$. The above integral consists of two parts: a radial part and an angular part which is an integration over the unit sphere. If we consider e.g. the case $s=2a$, $t=2b$, the angular integral simplifies to
\begin{displaymath}
 \int_{S^{m-1}} \overline{M_k(\underline{\omega})} \ \overline{\underline{\omega}^{2a}} \ \underline{\omega}^{2b} \ M_{\ell}(\underline{\omega}) \ w_k(\underline{\omega}) \ dS(\underline{\omega}) =  (-1)^{a+b} \ \int_{S^{m-1}} \overline{M_k(\underline{\omega})} \ M_{\ell}(\underline{\omega}) \ w_k(\underline{\omega}) \ dS(\underline{\omega})  .
\end{displaymath}
The case where $s$ and $t$ are both odd, i.e. $s=2a+1$ and $t=2b+1$, yields the same result as above.\\
On the other hand, if $s$ and $t$ have different parity, e.g. $s=2a$ and $t=2b+1$, the integral over the unit sphere vanishes, since outer and inner Dunkl monogenics are orthogonal (see Theorem \ref{orthomon}).

In what follows, we will restrict ourselves to spaces of polynomials of the type
\begin{displaymath}
R(M_k) = \left\lbrace p_n(\ux) \ M_k(\ux) = \sum_{j=0}^n a_j \ \ux^j \ M_k(\ux) \ | \ n \in \mathbb{N} \ , \ a_j \in \mathbb{R} \right\rbrace
\end{displaymath}
where $M_k$ is a left solid inner Dunkl monogenic of degree $k$, fixed once and for all, which satisfies
\begin{displaymath}
\int_{S^{m-1}} \overline{M_k(\underline{\omega})} \ M_k(\underline{\omega}) \ w_k(\underline{\omega}) \ dS(\underline{\omega}) = 1  .
\end{displaymath}

Inspired by the previous calculations and using the analytic continuation of the Gamma function, we are led to the following definition of a bilinear form on $R(M_k)$.
\begin{definition}~\\
The bilinear form $\langle.,.\rangle_{\alpha}$ (parameterized by $\alpha$) on $R(M_k)$ is defined by linear extension of
\begin{eqnarray*}
\langle \ux^{2s} M_k , \ux^{2t} M_k \rangle_{\alpha} & = & \frac{(-1)^{s+t}}{2} \ B \left( s+t+k+\frac{\mu}{2} , -\left( s+t+k+\frac{\mu}{2} \right) - \alpha \right)\\
\langle \ux^{2s+1} M_k , \ux^{2t} M_k \rangle_{\alpha} & = & 0\\
\langle \ux^{2s} M_k , \ux^{2t+1} M_k \rangle_{\alpha} & = & 0\\
\langle \ux^{2s+1} M_k , \ux^{2t+1} M_k \rangle_{\alpha} & = & \frac{(-1)^{s+t}}{2} \ B \left( s+t+k+\frac{\mu}{2}+1 , -\left( s+t+k+\frac{\mu}{2}+1 \right) - \alpha \right)  .
\end{eqnarray*}
\end{definition}
Note that this bilinear form is symmetric, but not always positive definite. Moreover, this bilinear form is well-defined if and only if $\alpha \notin \mathbb{N}$ and $\frac{\mu}{2} + \alpha \notin \pm \mathbb{N}$, due to the singularities $z=-n$, $n \in \mathbb{N}$ of the Gamma function $\Gamma(z)$.\\
Note that due to Lemma \ref{belangrijk} we have that $D_{\alpha} \lbrack R (M_k) \rbrack \subset R(M_k)$.
\begin{proposition}~\\
The operators $D_k$ and $D_{\alpha}$ are dual with respect to $\langle . \ , \ .\rangle_{\alpha}$, i.e.
\begin{displaymath}
\langle D_{\alpha} \lbrack p_i \ M_k \rbrack , p_j \ M_k \rangle_{\alpha} \ = \ \langle p_i \ M_k , D_k \lbrack p_j \ M_k \rbrack \rangle_{\alpha +1}
\end{displaymath}
with $p_i \ M_k$, $p_j \ M_k \in R(M_k)$.
\end{proposition}
\textit{Proof.} It is sufficient to prove the result for $\langle D_{\alpha} \lbrack \ux^{2s+1} M_k \rbrack , \ux^{2t} \ M_k \rangle_{\alpha}$, $\langle D_{\alpha} \lbrack \ux^{2s} M_k \rbrack , \ux^{2t+1} \ M_k \rangle_{\alpha}$, $\langle D_{\alpha} \lbrack \ux^{2s+1} M_k \rbrack , \ux^{2t+1} \ M_k \rangle_{\alpha}$ and $ \langle D_{\alpha} \lbrack \ux^{2s} M_k \rbrack , \ux^{2t} \ M_k \rangle_{\alpha}$.\\
By means of the definitions of the operator $D_{\alpha}$ and the bilinear form  $\langle .,. \rangle_{\alpha}$, we have consecutively
\begin{eqnarray*}
& & \hspace{-0,8cm} \langle D_{\alpha} \lbrack \ux^{2s+1} M_k \rbrack , \ux^{2t} \ M_k \rangle_{\alpha}\\
&= & -(2s+2k+\mu) \ \langle \ux^{2s} M_k , \ux^{2t} M_k \rangle_{\alpha} + (2s+2k+\mu+2 \alpha +2) \ \langle \ux^{2s+2} M_k , \ux^{2t} M_k \rangle_{\alpha}\\
& = & (2s+2k+\mu) \ \frac{(-1)^{s+t+1}}{2} \ \frac{\Gamma \left( s+t+k+\frac{\mu}{2} \right) \ \Gamma \left( - \left( s+t+k+\frac{\mu}{2} +\alpha \right) \right)}{\Gamma (- \alpha)}\\
& & + (2s+2k+\mu+2 \alpha+2) \frac{(-1)^{s+t+1}}{2} \ \frac{ \Gamma \left( s+1+t+k+\frac{\mu}{2} \right) \ \Gamma \left( - \left( s+t+k+\frac{\mu}{2} +\alpha +1 \right) \right)}{\Gamma (- \alpha)}.
\end{eqnarray*}
As $\Gamma(z+1) = z \ \Gamma(z)$, we can further simplify the above result
\begin{eqnarray*}
& & \hspace{-1,8cm} \langle D_{\alpha} \lbrack \ux^{2s+1} M_k \rbrack , \ux^{2t} \ M_k \rangle_{\alpha}\\
& = & -2 (\alpha +1) t \ \frac{(-1)^{s+t}}{2} \ \frac{\Gamma \left( s+t+k+\frac{\mu}{2} \right) \ \Gamma \left( - \left( s+t+k+\frac{\mu}{2} +\alpha +1 \right) \right)}{\Gamma (- \alpha)}\\
& = & 2t \ \frac{(-1)^{s+t}}{2} \ \frac{\Gamma \left( s+t+k+\frac{\mu}{2} \right) \ \Gamma \left( - \left(s+t+k+\frac{\mu}{2} + \alpha +1 \right) \right)}{\Gamma (-\alpha - 1)}\\
& = & \langle \ux^{2s+1} M_k , D_k \lbrack \ux^{2t} \ M_k \rbrack \rangle_{\alpha+1}  .
\end{eqnarray*}
The other three cases are treated similarly. \ \ \  \ \ \ $\square$

Now we come to the orthogonality relation of the Clifford-Gegenbauer polynomials.
\begin{theorem}~\\
If $s \not= t$, then
\begin{displaymath}
\langle {\mathcal G}_{s,\mu}^{\alpha}(M_k), {\mathcal G}_{t,\mu}^{\alpha}(M_k) \rangle_{\alpha} \ = \ 0  .
\end{displaymath}
\end{theorem}
\textit{Proof.} Suppose that $s>t$. The case where $s<t$ is similar. By means of the above proposition and Lemma \ref{belangrijk} we obtain consecutively
\begin{eqnarray*}
\langle {\mathcal G}_{s,\mu}^{\alpha}(M_k), {\mathcal G}_{t,\mu}^{\alpha}(M_k) \rangle_{\alpha} & = & \langle D_{\alpha} D_{\alpha +1} \ldots D_{\alpha +s-1} \lbrack M_k \rbrack \ , \ {\mathcal G}_{t,\mu}^{\alpha}(M_k) \rangle_{\alpha}\\
& = & \langle D_{\alpha +1} \ldots D_{\alpha +s-1} \lbrack M_k \rbrack \ , \ D_k \lbrack {\mathcal G}_{t,\mu}^{\alpha}(M_k) \rbrack \rangle_{\alpha+1}\\
& = & ...\\
& = & \langle M_k \ , \ D_k^s \lbrack {\mathcal G}_{t,\mu}^{\alpha}(M_k) \rbrack \rangle_{\alpha+s} \ = \  0  . \ \ \ \ \ \ \ \ \ \square
\end{eqnarray*}

\end{document}